\documentclass[12pt]{article}

\usepackage{amssymb}
\font\smallit=cmti10
\font\smalltt=cmtt10

\begin{document}
\begin{center}
{\uppercase{\bf Some Two Color, Four Variable Rado Numbers}}

\vskip 20pt

{\bf Aaron Robertson}\\
{\smallit Department of Mathematics,
Colgate University,
Hamilton, NY 13346}\\
{\smalltt aaron@math.colgate.edu}
\vskip 0pt
and
\vskip -5pt
{\bf Kellen Myers\footnote{
This work was done as part of a summer REU,
funded by Colgate University,
while the second author was an undergraduate at Colgate University,
under the directorship of the first author.
}}
\end{center}
\vskip 30pt
\begin{abstract}
There exists a minimum integer $N$ such
that any $2$-coloring of
$\{1,2,\dots,N\}$ admits a monochromatic solution to $x+y+kz =\ell w$ for
$k,\ell \in \mathbb{Z}^+$, where
$N$ depends on $k$ and $\ell$.  We determine
$N$ when $\ell-k \in \{0,1,2,3,4,5\}$,
for all $k,\ell$ for which
$\frac{1}{2}((\ell-k)^2-2)(\ell-k+1)\leq k \leq \ell-4$,
as well as for arbitrary $k$ when $\ell=2$.

\end{abstract}
\vskip 30pt
\section*{\normalsize 1. Introduction}

For $r \geq 2$, an $r$-coloring of the positive integers $\mathbb{Z}^+$ is an assignment 
$\chi : \mathbb{Z}^+ \rightarrow \{0,1,\dots,r-1\}$. Given a diophantine equation $\mathcal{E}$
in the variables $x_1,\dots,x_n$, we say a solution
$\{\bar{x}_i\}_{i=1}^n$ is monochromatic if $\chi(\bar{x}_i)=\chi(\bar{x}_j)$ for every $i,j$ pair. A
well-known theorem of Rado states that
a linear
homogeneous equation
$c_1x_1+\dots+c_nx_n=0$ with each
$c_i\in
\mathbb{Z}$ admits a monochromatic solution in $\mathbb{Z}^+$
 under any $r$-coloring of $\mathbb{Z}^+$,  for any $r \geq 2$,  if and
only if some nonempty subset of $\{c_i\}_{i=1}^n$ sums to zero.
Such an equation is said to satisfy Rado's regularity
condition. The smallest
$N$ such that any $r$-coloring of $\{1,2,\dots,N\}=[1,N]$ satisfies this condition is called the 
$r$-color Rado number
for the equation
$\mathcal{E}$. Rado also proved the following,
much lesser known, result.

\noindent{\bf Theorem 1} ({\it Rado [R]}) Let $\mathcal{E}=0$
 be a linear homogeneous equation
with integer coefficients. Assume that
$\mathcal{E}$ has at least 3 variables
with both positive and negative coefficients.
Then any $2$-coloring of $\mathbb{Z}^+$ admits
a monochromatic solution to $\mathcal{E}=0$.

For Rado's original proof (in German) see [R]; for a proof in English see [MR].

In this article we study the equation $x+y+kz =\ell w$ for positive integers
$k$ and $\ell$.  As such, we make the following notation.

\noindent
{\bf Notation} For $k$ a positive integer and
$j>-k$ an integer, let $\mathcal{E}(k,j)$ represent the equation
$$
x+y+kz=(k+j)w.
$$

\section*{\normalsize 2. A General Upper Bound}

\noindent
{\bf Definition} 
Let $\mathcal{E}$ be any equation that satisfies
the conditions in Theorem 1.  Denote by
$RR(\mathcal{E})$
the minimum integer $N$ such that any $2$-coloring of
$[1,N]$ admits  a monochromatic 
solution to $\mathcal{E}$.

Part of the following result is essentially a result due
to Burr and Loo [BL] who show that, for $j \geq 4$,
we have $RR(x+y=jw) = {j+1 \choose 2}$.
Their result was never published.  Below we present
our (independently derived) proof.

\noindent
{\bf Theorem 2 } Let $k,j\in \mathbb{Z}^+$ with $j \geq 4$. Then
$
RR(\mathcal{E}(k,j)) \leq
{j+1 \choose 2}.
$  
Furthermore,  for all $k \geq \frac{(j^2-2)(j+1)}{2}$,
we have
$
RR(\mathcal{E}(k,j)) =
{j+1 \choose 2}.
$ 

\noindent
{\it Proof.}
Let $\mathcal{F}$ denote the equation $x+y=jw$.
We will show that $
RR(\mathcal{F}) \leq
{j+1 \choose 2}.
$ 
Since any solution to $\mathcal{F}$ is also a
solution to $\mathcal{E}(k,j)$ for any $k \in \mathbb{Z}^+$,
the first statement will follow.

Assume, for a contradiction, that there
exists a $2$-coloring of $\left[1,{j+1 \choose 2} \right]$ with
no monochromatic solution to $\mathcal{F}$.
Using the colors red and blue,  
we let $R$ be the set of red integers and $B$ be the
set of blue integers.
We denote solutions of $\mathcal{F}$ by $(x,y,w)$
where $x,y,w \in \mathbb{Z}^+$.

Since $(x,y,w)=(j-1,1,1)$ solves $\mathcal{F}$,
we may assume that $1 \in R$ and $j-1 \in B$.  
 We separate
the proof into two cases.

\noindent
{\tt Case 1.} $j+1 \in B$. Assume $i \geq 1$ is red.
Considering $(1,ij-1,i)$ gives $ij-1 \in B$.
If $i \leq \left\lfloor \frac{j+1}{2}\right\rfloor$, this, in turn, gives
us
$i+1
\in R$ by considering
$(ij-1,j+1,i+1)$.
Hence $1,2,\dots,\left\lfloor\frac{j+1}{2} \right\rfloor+1$ are all red.
But then $\left(\left\lfloor \frac{j}{2}\right\rfloor,
\left\lfloor \frac{j+1}{2}\right\rfloor,1\right)$ is a red solution,
a contradiction.

\noindent
{\tt Case 2.} $j+1 \in R$. This implies that $j \left(\frac{j+1}{2}\right) \in B$.
Note also that the solutions $(1,j-1,1)$ and $\left(j\left(\frac{j-1}{2}\right),
j\left(\frac{j-1}{2}\right),j-1\right)$ give us $j-1 \in B$ and $j\left(\frac{j-1}{2}\right) \in R$.

First consider the case when $j$ is even.
 By considering $\left(\frac{j}{2},\frac{j}{2},1\right)$ we see
that $\frac{j}{2} \in B$.  Assume, for $i \geq 1$, that $\frac{(2i-1)j}{2} \in B$.
Considering $\left(j \left(\frac{j+1}{2}\right),\frac{(2i-1)j}{2}, \frac{j}{2}+i\right)$
we have $\frac{j}{2}+i \in R$.  This, in turn, implies that
$\frac{(2i+1)j}{2}\in B$ by considering
$\left(j \left(\frac{j-1}{2}\right),\frac{(2i+1)j}{2}, \frac{j}{2}+i\right)$.
Hence we have $\frac{j}{2}+i$ is red for
$1 \leq i \leq \frac{j}{2}$.  This gives us that $j-1 \in R$ (when $i=\frac{j}{2}-1$),
a contradiction.

Now consider the case when $j$ is odd.  We consider two subcases.

\noindent
{\tt Subcase i.} $j \in B$.  For $i \geq 1$, assume that $ij \in B$.
We obtain $\frac{j+1}{2} + i \in R$ by considering the solution
$\left(j\left(\frac{j+1}{2}\right),ij,\frac{j+1}{2} + i \right)$.
This gives us $(i+1)j \in B$ by considering
$\left(j\left(\frac{j-1}{2}\right),(i+1)j,\frac{j+1}{2} + i \right)$
Hence, we have that $j,2j,\dots,\left(\frac{j+1}{2}\right)j$ are all
blue, contradicting the deduction that $j \left(\frac{j-1}{2}\right) \in R$.

\noindent
{\tt Subcase ii.} $j \in R$.  We easily have $2 \in B$.
Next, we conclude that $j\left(\frac{j-3}{2}\right) \in R$
by considering $\left(j\left(\frac{j+1}{2}\right),j\left(\frac{j-3}{2}\right),j-1\right)$.
Then, the solution
$\left(j\left(\frac{j-3}{2}\right),j\left(\frac{j-1}{2}\right),j-2\right)$
gives us $j-2 \in B$.  
We use $\left(j\left(\frac{j-1}{2}\right),j,\frac{j+1}{2}\right)$ to see that
$\frac{j+1}{2}\in B$.  From $\left(j
\left(\frac{j+1}{2}\right)-2,2,\frac{j+1}{2}\right)$ we have $j\left(\frac{j+1}{2}\right)-2 \in R$.
To avoid $\left(j\left(\frac{j-1}{2}\right)+2,j\left(\frac{j+1}{2}\right)-2,j\right)$
being a red solution, we have $j\left(\frac{j-1}{2}\right)+2 \in B$.  This gives
us a contradiction; the solution $\left(j\left(\frac{j-1}{2}\right)+2,j-2,\frac{j+1}{2}\right)$
is blue.

This completes the proof of the first statement of the theorem.

For the proof of the second statement of the theorem,
we need only provide a lower bound of ${j+1 \choose 2} - 1$.
We first show that any solution to $x+y+kz = (k+j)w$ with $x,y,z,w < {j+1 \choose 2}$ must have $z=w$ when
$k \geq \frac{(j^2-2)(j+1)}{2}$. Assume, for a contradiction, that $z \neq w$.  If $z<w$, then
$(k+j)w \geq (k+j)(z+1) > kz + k$.  However,
$x+y < j(j+1)$ while $k > j(j+1)$ for $j \geq 3$.  Hence, $z \not < w$.
If $z>w$, then $(k+j)w \leq k(z-1) + j\left({j+1 \choose 2} - 1\right)$.
Since we have $x+y+kz = (k+j)w$ we now have
$2+kz \leq x+y+kz \leq k(z-1) + j\left({j+1 \choose 2} - 1\right)$.
Hence, $k \leq j\left({j+1 \choose 2} - 1\right)-2$, contradicting
the given bound on $k$.  Thus, $z=w$.

Now, any solution to $x+y+kz = (k+j)w$ with $z=w$ is a solution
to $x+y=jw$.  From Burr and Loo's result, there exists a
$2$-coloring of $\left[1,{j+1 \choose 2}-1\right]$ with no
monochromatic solution to $x+y=jw$.  This provides us with a
$2$-coloring with no monochromatic solution to $x+y+kz=(k+j)w$, thereby
finishing the proof of the second statement.
 \hfill $\Box$

\section*{\normalsize 3.  Some Specific Numbers}

In this section we determine the exact
values for $RR(\mathcal{E}(k,j))$ for $j \in \{0,1,2,3,4,5\}$, most of which are cases
not covered by Theorem 2.

\noindent
{\bf Theorem 3} For $k\geq 2$,
$$
RR(\mathcal{E}(2k,0)) = 2k \mbox{ and } RR(\mathcal{E}(2k-1,0)) = 3k-1.$$
Furthermore $RR(\mathcal{E}(2,0))=5$ and $RR(\mathcal{E}(1,0))=11$. 

\noindent
{\it Proof.} The cases $RR(\mathcal{E}(2,0))=5$, $RR(\mathcal{E}(4,0))=4$,
and $RR(\mathcal{E}(3,0)=5$ are easy calculations, as is
$RR(\mathcal{E}(1,0))=11$, which first appeared in [BB].  Hence,
we may assume $k \geq 3$ in the following arguments.

We start with $RR(\mathcal{E}(2k,0)) = 2k$.  To show that
$RR(\mathcal{E}(2k,0)) \geq 2k$ consider the $2$-coloring
of $[1,2k-1]$ defined by coloring the odd integers red and
the even integers blue.  To see that there is no monochromatic
solution to $x+y+2kz = 2kw$, note that we must have
$2k \mid (x+y)$.  This implies that $x+y=2k$ since $x,y \leq 2k-1$.
Thus, $w=z+1$.  However, no 2 consecutive integers have the same
color.  Hence, any solution to $\mathcal{E}(2k,0)$ is necessarily
bichromatic.

Next, we show that $RR(\mathcal{E}(2k,0)) \leq 2k$.
Assume, for a contradiction, that there exists a $2$-coloring
of $[1,2k]$ with no monochromatic solution to our equation.
Using the colors red and blue, we may assume, without loss of generality, that
$k$ is red.  This gives us that $k-1$ and $k+1$ are blue, by considering
$(x,y,z,w)=(k,k,k-1,k)$ and $(k,k,k,k+1)$.  Using these in the
solution $(2k,2k,k-1,k+1)$ we see that $2k$ must be red, which
implies that $k-2$ is blue (using $(2k,2k,k-2,k)$).
However, this gives the blue solution $(k-1,k+1,k-2,k-1)$, a contradiction.

We move on to $RR(\mathcal{E}(2k-1,0))$.  To show
that $RR(\mathcal{E}(2k-1,0)) \leq 3k-1$ consider the
following $2$-colorings of $[1,3k-2]$, dependent on $k$
(we use $r$ and $b$ for red and blue, respectively):
$$
\begin{array}{ll}
brrbbrrbb \dots bbr &\mbox{ if }k \equiv 0 \,(\bmod \, 4)\\
rrbbrrbb \dots bbr &\mbox{ if }k \equiv 1 \,(\bmod \, 4)\\
brrbbrrbb \dots rrb  &\mbox{ if }k \equiv 2 \,(\bmod \, 4)\\
rrbbrrbb \dots rrb  &\mbox{ if }k \equiv 3 \,(\bmod \, 4).
\end{array}
$$
Since we need $(2k-1) \mid (x+y)$ and $x,y \leq 3k-2$, we
have $x+y \in \{2k-1,4k-2\}$.  By construction, if
$x+y=2k-1$, then $x$ and $y$ have different colors.  Hence,
the only possibility is $x+y=4k-2$.  But then
$w=z+2$ and we see that $w$ and $z$ must have different colors.

Next, we show that $RR(\mathcal{E}(2k-1,0)) \leq 3k-1$.
Assume, for a contradiction, that there exists a $2$-coloring
of $[1,3k-1]$ with no monochromatic solution to our equation.
Using the colors red and blue, we may assume, without loss of generality, that
$2k-1$ is red.  To avoid $(2k-1,2k-1,z,z+2)$ being a red solution,
we see that $2k+1$ and $2k-3$ are blue
(using $z=2k-1$ and $2k-3$, respectively).

If $2k$ is red, then $2k-2$ is blue (using $(2k-2,2k,2k-2,2k)$).
From $(3k-1,3k-1,2k-2,2k+1)$ we see that $3k-1$ is red.  This,
in turn, implies that $2k-4$ is blue (using $(3k-1,3k-1,2k-4,2k-1)$).
So that $(2k-4,2k+2,2k-4,2k-2)$ is not a blue solution, we require
$2k+2$ to be red.  But then $(2k-1,2k-1,2k,2k+2)$ is a red
solution, a contradiction.

If $2k$ is blue, then $2k-2$ must be red. 
So that $(3k-1,3k-1,2k-3,2k)$ is not a blue solution,
we have that $3k-1$ is red.  Also, $2k+2$ must be
red by considering $(2k-3,2k+1,2k,2k+2)$.  But this
implies that $(3k-1,3k-1,2k-1,2k+2)$ is a red solution,
a contradiction.
\hfill $\Box$

We proceed with a series of results for the cases
$j=1,3,4,5$.  When $j=2$, the corresponding number
is trivially $1$ for all $k \in \mathbb{Z}^+$.

Below, we will call a coloring of $[1,n]$ {\it valid}
if it does not contain a monochromatic
solution to $\mathcal{E}(k,j)$.

\noindent
{\bf Theorem 4}  For $k \in \mathbb{Z}^+$,
$$
RR(\mathcal{E}(k,1)) =
\left\{
\begin{array}{ll}
4&\mbox{ for }k\leq 3\\
5&\mbox{ for }k \geq 4.
\end{array}
\right.
$$

\noindent
{\it Proof.} Assume, for a contradiction,
that there exists a $2$-coloring of $[1,5]$
with no monochromatic solution to $x+y+kz=(k+1)w$.
We may assume that $1$ is red.
Considering the solutions $(1,1,2,2), (2,2,4,4), (1,3,4,4),$
and $(2,3,5,5)$, in order, we find that
$2$ is blue, $4$ is red, $3$ is blue, and $5$
is red.  But then $(1,4,5,5)$ is a red solution,
a contradiction.  Hence, $RR(\mathcal{E}(k,1)) \leq 5$
for all $k \in \mathbb{Z}^+$.

We see from the above argument that the only
valid colorings of $[1,3]$ (assuming,
without loss of generality, that $1$ is red) are
$rbr$ and $rbb$ (where we use $r$ for red and
$b$ for blue).  Furthermore, the only valid
coloring of $[1,4]$ is $rbbr$.  We use these colorings to
finish the proof.

First consider the valid coloring $rbr$.
The possible values of $x+y+kz$ when $x,y,z$ are
all red form the set $\{k+2,k+4,k+6,3k+2,3k+4,3k+6\}$.
The possible values when $x,y,z$ are all blue is $2k+4$.
The possible values of $(k+1)w$ when $w$ is red form the
set $\{k+1,3k+3\}$; when $w$ is blue, $2k+2$ is the
only possible value.  We denote these results by:
\renewcommand{\arraystretch}{1.25}
$$
\begin{array}{rl}
R_{x,y,z} \!\!&= \{k+2,k+4,k+6,3k+2,3k+4,3k+6\}\\
B_{x,y,z} \!\!&= \{2k+4\}\\
R_w \!\!&= \{k+1,3k+3\}\\
B_w \!\!&=\{2k+2\}.
\end{array}
$$

Next, we determine those values of $k$, if any,
for which $R_{x,y,z} \cap R_w \neq \emptyset$
or $B_{x,y,z} \cap B_w \neq \emptyset$.
Clearly, there is no such $k$ for these sets.
Hence, we conclude that $RR(\mathcal{E}(k,1)) \geq 4$
for all $k$.  (We need not consider the valid coloring
$rbb$ since we now know that 
$RR(\mathcal{E}(k,1)) \geq 4$
for all $k$.)

We  move on to the valid coloring of $[1,4]$, which
is $rbbr$.  We find that
$$
\begin{array}{rl}
R_{x,y,z} \!\!&= \{k+2,k+5,k+8,4k+2,4k+5,4k+8\}\\
B_{x,y,z} \!\!&= \{2k+4,2k+5,2k+6,3k+4,3k+5,3k+6\}\\
R_w \!\!&= \{k+1,4k+4\}\\
B_w \!\!&=\{2k+2,3k+3\}.
\end{array}
$$
\renewcommand{\arraystretch}{1.0}

We see that $B_{x,y,z} \cap B_w \neq \emptyset$
when $k=1$ ($2k+4=3k+3$), $k=2$ ($2k+5=3k+3$),
and $k=3$ ($2k+6=3k+3$).  For all other
values of $k$, $B_{x,y,z} \cap B_w = \emptyset$
and $R_{x,y,z} \cap R_w = \emptyset$.
Hence, we conclude that $RR(\mathcal{E}(k,1)) \geq 5$
for $k \geq 3$, while, since
$rbbr$ is the only valid coloring of $[1,4]$,
$RR(\mathcal{E}(k,1)) \leq 4$ for $k=1,2,3$.
This completes the proof of the theorem.
\hfill{$\Box$}

The proofs below refer to the small
Maple package {\tt FVR}.  The description of
{\tt FVR} is given in Section 3.1, which follows the next
3 theorems.

\noindent
{\bf Theorem 5}  For $k \in \mathbb{Z}^+$,
$$
RR(\mathcal{E}(k,3)) =
\left\{
\begin{array}{ll}
4&\mbox{ for }k \leq 5 \mbox{ and } k=7\\
6&\mbox{ for }k =8,11\\
9&\mbox{ for }k=6,9,10 \mbox { and }k \geq 12
\end{array}
\right.
$$

\noindent
{\it Proof.}  The method of proof is the same as that
for Theorem 4, but we will work it out in
some detail commenting on the use of the Maple
package {\tt FVR}.

 It is easy to check that
the only valid $2$-colorings (using $r$ for red,
$b$ for blue, and assuming that $1$ is red) of $[1,n]$
for $n=4,5,\dots,8$ are as in the following table.
The determinations of $R_{x,y,z}, R_w, B_{x,y,z},$ and
$B_w$ are equally easy.
\small
\renewcommand{\arraystretch}{1.25}
$$
\begin{array}{r|l|l}
n&\mbox{coloring}&\mbox{sets}\\ \hline
4&rbrr& R_{x,y,z}=\{ik+j: i=1,3,4; j=2,4,5,6,7,8\}; R_w=\{i(k+3): i=1,3,4\}\\
&&B_{x,y,z}=\{2k+2\};B_w=\{2k+6\}\\
5&rbrrb& R_{x,y,z}=\{ik+j: i=1,3,4; j=2,4,5,6,7,8\}; R_w=\{i(k+3): i=1,3,4\}\\
&&B_{x,y,z}=\{ik+j: i=2,5;j=2,4,7,10\}; B_w=\{2k+6,5k+15\}\\
6&rbrrbb& R_{x,y,z}=\{ik+j: i=1,3,4; j=4,5,6,7,8\}; R_w=\{i(k+3): i=1,3,4\}\\
&&B_{x,y,z}=\{ik+j: i=2,5,6;j=4,7,8,10,11,12\}; B_w=\{2k+6,5k+15,6k+18\}\\
7&rbrrbbr& R_{x,y,z}=\{ik+j: i=1,3,4,7; j=4,\dots,8,10,11,14\}; R_w=\{i(k+3): i=1,3,4,7\}\\
&&B_{x,y,z}=\{ik+j: i=2,5,6;j=4,7,8,10,11,12\}; B_w=\{2k+6,5k+15,6k+18\}\\
8&rbrrbbrb& R_{x,y,z}=\{ik+j: i=1,3,4,7; j=4,\dots,8,10,11,14\}; R_w=\{i(k+3): i=1,3,4,7\}\\
&&B_{x,y,z}=\{ik+j: i=2,5,6;j=4,7,8,10,\dots,14,16\}; B_w=\{i(k+3): i=2,5,6,8\}
\end{array}
$$
\renewcommand{\arraystretch}{1.0}

\normalsize
The sets $R_{x,y,z}, R_w, B_{x,y,z},$ and
$B_w$ are automatically found by {\tt FVR}, which
then gives us the values of $k$ that
induce a nonempty intersection of
either $R_{x,y,z} \cap R_w$ or
$B_{x,y,z} \cap B_w$.  For completeness, we give the
details.

For the coloring $rbrr$, we have
$R_{x,y,z} \cap R_w \neq \emptyset$ when
$k=1$ ($k+3=3k+5$), $k=2$ ($3k+9=4k+7$),
$k=3$ ($3k+9=4k+6$), $k=4$ ($3k+9=4k+5$),
$k=5$ ($3k+9=4k+4$), and
$k=7$ ($3k+9=4k+2$).  Since $rbrr$ is the only
valid coloring of $[1,4]$, we have
$RR(\mathcal{E}(k,3))=4$ for $k=1,2,3,4,5,7$.

For the coloring $rbrrb$, we have no new additional
elements in $R_{x,y,z} \cap R_w$.  Hence, any possible
additional intersection point comes from 
$B_{x,y,z} \cap B_w$.  However,
$B_{x,y,z} \cap B_w = \emptyset$ for all $k\in \mathbb{Z}^+$.
Hence, $RR(\mathcal{E}(k,3)) \geq 6$ for $k \in \mathbb{Z}^+ \setminus \{1,2,3,4,5,7\}$.

For $rbrrbb$, we again have no new additional elements in
the red intersection.  We do, however, have
additional elements in $B_{x,y,z} \cap B_w$.
When $k=8$ ($5k+15=6k+7$) and $k=11$ ($5k+15=6k+4$) we
have a blue intersection.  We conclude that
$RR(\mathcal{E}(k,3))=6$ for $k=8,11$ and 
$RR(\mathcal{E}(k,3)) \geq 7$ for $k \in \mathbb{Z}^+ \setminus \{1,2,3,4,5,7,8,11\}$.

Considering $rbrrbbr$, we have no new additional elements in
$B_{x,y,z} \cap B_w$.  Furthermore, we have no new additional
intersection points in $R_{x,y,z} \cap R_w$.
Hence, $RR(\mathcal{E}(k,3)) \geq 8$ for $k \in \mathbb{Z}^+ \setminus \{1,2,3,4,5,7,8,11\}$.

Lastly, we consider $rbrrbbrb$, which gives no new additional
elements in $R_{x,y,z} \cap R_w$.
Furthermore, we have no new additional intersection points in
$B_{x,y,z} \cap B_w$.
Thus, $RR(\mathcal{E}(k,3)) \geq 9$ for $k \in \mathbb{Z}^+ \setminus \{1,2,3,4,5,7,8,11\}$.

Analyzing the valid coloring of $[1,8]$
we see that we cannot extend it to a valid
coloring of $[1,9]$.  Hence,
$RR(\mathcal{E}(k,3)) \leq 9$ for all $k$ so
that 
$RR(\mathcal{E}(k,3))= 9$ for $k \in \mathbb{Z}^+ \setminus \{1,2,3,4,5,7,8,11\}$.
\hfill{$\Box$}

\noindent
{\bf Theorem 6}  For $k \in \mathbb{Z}^+$,
$$
RR(\mathcal{E}(k,4)) =
\left\{
\begin{array}{rl}
3&\mbox{ for }k =2,3,4 \\
5&\mbox{ for }k =6,7,8,10,11,14\\
6&\mbox{ for }k= 5,9,12,13,15,18\\
8&\mbox{ for }k=17,19,22\\
9&\mbox{ for }k=1,23,24\\
10&\mbox{ for }k=16,20,21 \mbox{ and } k\geq 25.
\end{array}
\right.
$$

\noindent
{\it Proof.}  Use the Maple package {\tt FVR} with the following valid
colorings (which are easily obtained):

$$
\begin{array}{r|l}
n&\mbox{valid colorings}\\ \hline
3&rbb\\
4&rbbr\\
5&rbbrr\\
6&rbbrrr\\
7&rbbrrrr, rbbrrrb\\
8&rbbrrrrb, rbbrrrbb\\
9&rbbrrrrbr, rbbrrrbbr\\
10& \mbox{none}
\end{array}
$$

Note that if $[1,n]$ has more than one valid coloring, we can 
conclude that
$RR(\mathcal{E}(k,4)) \leq n$ for $k=\hat{k}$ only if
$\hat{k}$ is an intersection point for {\it all} valid colorings.
Otherwise, there exists a coloring of $[1,n]$ that avoids
monochromatic solutions to $\mathcal{E}(k,4)$ when $k=\hat{k}$.
\hfill{$\Box$}

\noindent
{\bf Theorem 7}  For $k \in \mathbb{Z}^+$,
$$
RR(\mathcal{E}(k,5)) =
\left\{
\begin{array}{rl}
4&\mbox{ for }k =1,2,3 \\
6&\mbox{ for }k =4,13,14\\
7&\mbox{ for }k= 16,17,18,23\\
8&\mbox{ for }5 \leq k \leq 12 \mbox{ and } k=21\\
10&\mbox{ for }k=19,24,26,27,28,29,33\\
11&\mbox{ for }k=22,30,31,32,34,36,37,38,39,41,42,43,48\\
12&\mbox{ for }k=15,35,44,46,47,53\\
13&\mbox{ for }k=51,52\\
15&\mbox{ for }k=20,25,40,45,49,50 \mbox{ and } k\geq 54.
\end{array}
\right.
$$

\noindent
{\it Proof.}   Use the Maple package {\tt FVR} with the following valid
colorings (which are easily obtained):

$$
\begin{array}{r|l}
n&\mbox{valid colorings}\\ \hline
4&rrrb, rrbb, rbrb, rbbb\\
5&rrbbb,rbrbr, rbbbr\\
6&rrbbbr, rbrbrr, rbbbrr\\
7&rrbbbrr, rrbbbrb, rbrbrrr, rbrbrrb, rbbbrrr\\
8&rrbbbrrb, rbrbrrbr, rbbbrrrr\\
9&rrbbbrrbb, rbrbrrbrb, rbbbrrrrr. rbbbrrrrb\\
10&rrbbbrrbbr, rbbbrrrrbr, rbbbrrrrrr\\
11&rrbbbrrbbrr, rbbbrrrrrrr, rbbbrrrrbrr\\
12&rrbbbrrbbrrr, rbbbrrrrrrrr, rbbbrrrrbrrr\\
13&rrbbbrrbbrrrr, rrbbbrrbbrrrb\\
14&rrbbbrrbbrrrrr,rrbbbrrbbrrrrb, rrbbbrrbbrrrbr,rrbbbrrbbrrrbb\\
15& \mbox{none}.
\end{array}
$$
\hfill{$\Box$}

\subsection*{\normalsize 3.1 About {\tt FVR}}

In the above theorems, we find
our lower bounds by considering all valid
colorings of $[1,n]$ for some $n\in \mathbb{Z}^+$
and deducing the possible elements that
$x+y+kz$ can be when $x,y$, and $z$ are monochromatic
and the possible elements that $(k+j)w$ can be,
i.e., determining $R_{x,y,z}, B_{x,y,z}, R_w,$ and $B_w$.
We then looked for intersections that would
make $(x,y,z,w)$ a monochromatic solution.
The intersections are specific values of $k$ which show
that the given coloring has monochromatic solutions
for these values of $k$.

This process has been automated
in the Maple package {\tt FVR}, which is
available from the first author's
website\footnote{\tt http://math.colgate.edu/$\sim$aaron/programs.html}.
The input is a list of all valid colorings of
$[1,n]$.  The output is a list of values of
$k$ for which we have monochromatic solutions.
By increasing $n$ we are able to determine
the exact Rado numbers for all $k \in \mathbb{Z}^+$.
An example of this is explained in detail in
the proof of the next theorem.

\section*{\normalsize 4. A Formula for $x+y+kz=2w$}

In [HS] and [GS] a formula for, in particular, $x+y+kz=w$ is given:
$RR(x+y+kz=w) = (k+1)(k+4)+1$.  In this section we provide
a formula for the next important equation of this form,
namely the one in this section's title.  To the best of
our knowledge this is the first formula given for
a linear homogeneous equation $\mathcal{E}=0$ 
 of more
than three variables with a negative coefficient
not equal to $-1$ (assuming, without loss of generality,
 at least as many positive
coefficients as negative ones)
that does not satisfy Rado's regularity condition.

\noindent
{\bf Theorem 8}  For $k \in \mathbb{Z}^+$,
$$
\renewcommand{\arraystretch}{1.5}
RR(x+y+kz=2w) = \left\{
\begin{array}{ll}
\frac{k(k+4)}{4}+1& \mbox{if } k \equiv 0 \,(\mbox{mod }4)\\
\frac{(k+2)(k+3)}{4}+1& \mbox{if } k \equiv 1 \,(\mbox{mod }4)\\
\frac{(k+2)^2}{4}+1& \mbox{if } k \equiv 2 \,(\mbox{mod }4)\\
\frac{(k+1)(k+4)}{4}+1& \mbox{if } k \equiv 3 \,(\mbox{mod }4).
\end{array}
\right.
$$

\noindent
{\it Proof.}  We begin with the lower bounds.  Let $N_i$ be one less
than the stated formula for $k \equiv i\, (\mbox{mod }4)$,
with $i \in \{0,1,2,3\}$.
We will provide $2$-colorings of $[1,N_i]$, for $i=0,1,2,3$, that admit
no monochromatic solution to $x+y+kz=2w$.

For $i=0$, color all elements in $\left[1,\frac{k}{2}\right]$ red and
all remaining elements blue.  If we assume $x,y,$ and $z$ are
all red, then $x+y+kz \geq k+2$ so that for any solution we have
$w \geq \frac{k}{2}+1$.  Thus, there is no red solution.
If we assume $x,y,$ and $z$ are all blue,
then $x+y+kz \geq \frac{k^2}{2}+2k+2 = 2\left(\frac{k(k+4)}{4}+1\right) >
2N_0$, showing that there is no blue solution.

For $i=1$, color $N_1$ and all elements in $\left[1,\frac{k+1}{2}\right]$ red.
Color the remaining elements blue.  Similarly to the last case,
we have no blue solution
since $x+y+kz > 2(N_1-1)$.  If we assume $x,y,z$ and $w$ are all
red, then we cannot have all of $x,y,z$ in
$\left[1,\frac{k+1}{2}\right]$. If we do, since $k$ is odd
(so that we must have $x+y$ odd), then
$x+y+kz \geq 1+2+k(1) = k +3 > 2 \left(\frac{k+1}{2}\right)$. Thus, $w$
must be blue. Now  we assume, without loss of generality,
that $x = N_1$.  In this situation, we must have
$w=N_1$.  Hence, since
$N_1+y+kz=2N_1$ we see that $y+kz = N_1$.  Hence, $y,z \leq
\frac{k+1}{2}$. But then $y+kz \leq \frac{k^2+2k+1}{2} <
N_1$, a contradiction. Hence, there is no red solution under this
coloring.

The cases $i=2$ and $i=3$ are similar to the above cases.
As such, we provide the colorings and leave the details to the reader.
For $i=2$, we color $N_2$ and all elements in $\left[1,\frac{k}{2}\right]$ red,
while the remaining elements are colored blue.
For $i=3$, color all elements in $\left[1,\frac{k+1}{2}\right]$ red and
all remaining elements blue. 

We now turn to the upper bounds.  
We let $M_i$ be equal to the stated formula for $k \equiv i\, (\mbox{mod }4)$,
with $i \in \{0,1,2,3\}$.  We employ a ``forcing" argument to
determine the color of certain elements.  We let $R$ denote the
set of red elements and $B$ the set of blue elements.
We denote by a $4$-tuple $(x,y,z,w)$ a solution to $x+y+kz=2w$.
In each of the following cases assume, for a contradiction, that there
exists a $2$-coloring of $[1,M_i]$ with no
monochromatic solution to the equation.
In each case we assume $1 \in R$.

\noindent
{\tt Case 1.} $k \equiv 0 \, (\mbox{mod }4)$. 
We will first show that $2 \in R$.  Assume, for a contradiction,
that $2 \in B$.  Then $2k+2 \in R$ by considering
$(2,2k+2,2,2k+2)$.  Also, $k+1 \in B$ comes
from the similar solution $(1,k+1,1,k+1)$.
Now, from $(3k+3,1,1,2k+2)$ we have $3k+3 \in B$.
As a consequence, we see that $3 \in R$ by considering
$(3,3k+3,3,3k+3)$.  From here we use $(3k+1,3,1,2k+2)$ to
see that $3k+1 \in B$.  But then
$(3k+1,k+1,2,3k+1)$ is a blue solution, a contradiction.
Hence, $2 \in R$.

Now, since $1,2 \in R$, in order for
 $\left(1,1,1,\frac{k}{2}+1\right)$ not to be monochromatic, we
have $\frac{k}{2}+1 \in B$.  Similarly, $\left(2,2,1,\frac{k}{2}+2\right)$
gives  $\frac{k}{2}+2 \in B$.  Consequently, so that
$\left(\frac{k}{2}+1,\frac{k}{2}+1,\frac{k}{2}+1,\frac{k^2}{4}+k+1\right)$
is not monochromatic,
we have $\frac{k^2}{4}+k+1 \in R$.

Our next goal is to show that $\frac{k}{2} \in R$.  So that
$\left(\frac{k^2+2k}{4}+1,\frac{k^2+2k}{4}+1,1,\frac{k^2}{4}+k+1 \right)$ is not red,
we have $\frac{k^2+2k}{4}+1 \in B$.  In turn, to avoid
$\left(\frac{k}{2}+1,\frac{k}{2}+1,\frac{k}{2},\frac{k^2+2k}{4}+1 \right)$
being blue, we have $\frac{k}{2} \in R$, as desired.

So that $\left(\frac{k}{2},\frac{k}{2},1,k\right)$
and $\left(\frac{k}{2},\frac{k}{2},\frac{k}{2},\frac{k^2+2k}{4}\right)$ are not red, we have 
$k,\frac{k^2+2k}{4} \in B$.
Using these in $\left(k,k,\frac{k}{2}-1,\frac{k^2+2k}{4}\right)$ gives 
$\frac{k}{2}-1 \in R$.
Since $\frac{k}{2}$ and $\frac{k}{2} -1$ are both red,
$\left(\frac{k}{2} -1,\frac{k}{2} -1,\frac{k}{2}-1, \frac{k^2}{4}-1\right)$
gives us $\frac{k^2}{4}-1 \in B$ while
$\left(\frac{k}{2},\frac{k}{2},\frac{k}{2}-1,\frac{k^2}{4}\right)$
gives us $\frac{k^2}{4} \in B$. 
This, in turn, gives us $\frac{k}{4} \in R$ by
considering $\left(\frac{k^2}{4},k,\frac{k}{4},\frac{k^2+2k}{4}\right)$.

Now, from
$\left(\frac{k^2+4k}{4}+1,1,\frac{k}{4}+1,\frac{k^2+4k}{4}+1\right)$
we have $\frac{k}{4}+1 \in B$.  We use this in the two solutions
$\left(\frac{k^2}{4}-1,\frac{k}{2}+1,\frac{k}{4}+1,\frac{k^2+3k}{4}\right)$
and
$\left(\frac{k^2}{4},k,\frac{k}{4}+1,\frac{k^2+4k}{4}\right)$
to find that $\frac{k^2+3k}{4},\frac{k^2+4k}{4} \in R$.
But this gives us the red solution
$\left(\frac{k^2+4k}{4},\frac{k}{2},\frac{k}{4},\frac{k^2+3k}{4}\right)$,
a contradiction.

\noindent
{\tt Case 2.} $k \equiv 1 \, (\mbox{mod }4)$.
The argument at the beginning of Case 1 holds for
this case, so we have $2 \in R$. 
We consider two subcases.

\noindent
{\tt Subcase i.} $ \frac{k+3}{4}\in R$. From
$\left(\frac{k^2+4k+3}{4}, \frac{k+3}{4}, \frac{k+3}{4},\frac{k^2+4k+3}{4}\right)$
we have $\frac{k^2+4k+3}{4} \in B$.  This gives
us $\frac{k+1}{2} \in R$ by considering
$\left(k+1,\frac{k+1}{2},\frac{k+1}{2},\frac{k^2+4k+3}{4}\right)$
(where $k+1 \in B$ comes from $(1,1,2,k+1)$).
We also have, from $(1,2,1,\frac{k+3}{2})$, that $\frac{k+3}{2} \in B$.
Consequently, $\frac{k^2+5k+6}{4} \in R$ so that
$\left(\frac{k+3}{2},\frac{k+3}{2},\frac{k+3}{2},\frac{k^2+5k+6}{4}\right)$
is not monochromatic.

We next have that $\frac{k+5}{2} \in B$ so that
$\left(3,2,1,\frac{k+5}{2}\right)$ is not monochromatic
(we may assume that $k \geq 9$).  Hence,
$\frac{k^2+5k+10}{4} \in R$ by considering
$\left(\frac{k+5}{2},\frac{k+5}{2},\frac{k+3}{2},\frac{k^2+5k+10}{4}\right)$.
But this gives us the monochromatic solution
$\left(\frac{k^2+5k+10}{4},\frac{k+1}{2},\frac{k+3}{4},\frac{k^2+5k+6}{4}\right)$,
a contradiction.

\noindent
{\tt Subcase i.} $ \frac{k+3}{4}\in B$.  Via arguments similar to those
in Subcase i, we have $\frac{k^2+4k+3}{4}, \frac{k^2+5k+6}{4} \in R$.
From
$\left(\frac{k^2+4k+3}{4},\frac{k^2+2k+9}{4},1,\frac{k^2+5k+6}{4}\right)$
we have $\frac{k^2+2k+9}{4} \in B$.  This gives us $\frac{k^2+15}{4} \in
R$ by considering $\left(\frac{k^2+15}{4},\frac{k+3}{4},\frac{k+3}{4},
\frac{k^2+2k+9}{4}\right)$.

We now show that $\frac{k-1}{4} \in R$ by showing
that for any $i \leq \frac{k-1}{4}$ we must
have $i \in R$.  To this end, assume, for a contradiction,
that $i-1 \in R$ but $i\in B$ (where $i \geq 3$).
From $(i,ik+i,i,ik+i)$ we have $ik+i \in R$.  In turn we
have $(i+1)k+i \in B$ by considering $(ik+i,ik+i,2,(i+1)k+i)$.
We next see from $((i+1)k+i,i,i+1,(i+1)k+i)$ that $i+1 \in R$.
Using our assumption that $i-1 \in R$ in $(i-1,i+1,2,k+i)$
we have $k+i \in B$.  But then
$((i+1)k+i,k+i,i,(i+1)k+i)$ is a blue solution,
provided $(i+1)k+i \leq M_1$, which by the bound
given on $i$ is valid.  By applying this argument
to $i=3,4,\dots, \frac{k-1}{4}$, in order, we see that all
positive integers less than or equal to $ \frac{k-1}{4}$
must be red.  In particular, $ \frac{k-1}{4} \in R$.

Using $\frac{k-1}{4} \in R$ in $\left(
\frac{k^2+15}{4},\frac{k+15}{4},\frac{k-1}{4},\frac{k^2+15}{4}\right)$
we have $\frac{k+15}{4} \in B$.  This, in turn, gives us
$\frac{k^2+4k+15}{4} \in R$ by considering
$\left(
\frac{k^2+4k+15}{4},\frac{k+15}{4},\frac{k+3}{4},\frac{k^2+4k+15}{4}\right)$.
For our contradiction, we see now that
$\left(
\frac{k^2+4k+15}{4},\frac{k^2+15}{4},1,\frac{k^2+4k+15}{4}\right)$
is a red solution.

\noindent
{\tt Case 3.} $k \equiv 2 \, (\mbox{mod }4)$.
From Case 1 we have $\frac{k}{2},\frac{k^2}{4}+k+1 \in R$ and
$\frac{k}{2}+1,\frac{k^2}{4}-1,\frac{k^2+2k}{4} \in B$.
From $\left(1,1,\frac{k}{2}+2,\frac{k^2}{4}+k+1\right)$ we see that
$\frac{k}{2}+2 \in B$.  This gives us $\frac{k^2}{4}+k + 2 \in R$
by considering $\left(\frac{k}{2}+2 ,\frac{k}{2}+2,\frac{k}{2}+1,
\frac{k^2}{4}+k + 2\right)$.  Using this fact in
$\left(\frac{k^2}{4}+k + 2,\frac{k}{2},\frac{k+2}{4}, \frac{k^2}{4}+k + 1\right)$
we have $\frac{k+2}{4} \in B$.  But then
$\left(\frac{k^2}{4}-1,\frac{k}{2}+1,\frac{k+2}{4},\frac{k^2+2k}{4}\right)$
is a blue solution, a contradiction.

\noindent
{\tt Case 4.} $k \equiv 3 \, (\mbox{mod }4)$. Let $i \in R$.
From Case 2 we may assume $i \geq 2$ so
that $1,2,\dots,i$ are all red,
By considering $(i,ik+i,i,ik+i)$ we have $ik+i \in
B$ so that we may assume
$k+1,2k+2,\dots,(i-1)k+(i-1),ik+i$ are all blue.
Since $(i+1,(i-1)k+i-1,i+1,ik+i)$ is a solution, we have
$i+1 \in R$.  Hence, $i \in R$ for $1 \leq i \leq \frac{k+5}{4}$.
In particular, $\frac{k+5}{4} \in R$.
From Case 2 we also have
$\frac{k+3}{2},\frac{k+5}{2}\in B$.  
By considering $\left(\frac{k+3}{2},\frac{k+5}{2},\frac{k+3}{2},\frac{k^2+5k+8}{4}\right)$
we have $\frac{k^2+5k+8}{4} \in R$.
But then
$\left(\frac{k^2+5k+8}{4},2,\frac{k+5}{4},\frac{k^2+5k+8}{4}\right)$
is a red solution, a contradiction.

\section*{\normalsize 4. Conclusion}

The next important numbers to determine are in the first row
of Tables 1.  As such,
it would be nice to have a formula for $RR(x+y+z=\ell w)$.
We have been unable to discover one.

By analyzing Table 1, given below, we have noticed,
to some extent, certain patterns that emerge.  In particular,
we make the following conjecture.

\noindent
{\bf Conjecture}  For $\ell \geq 2$ fixed and $k \geq \ell+2$,
we have 
$$RR(x+y+kz = \ell w) = \left(\left\lfloor
\frac{k+\ell+1}{\ell}\right\rfloor\right)^2 +
O\left(\frac{k}{\ell^2}\right),$$
where the ``$O\left(\frac{k}{\ell^2}\right)$ part" depends on
the residue class of $k$ modulo $\ell^2$.

For a concrete conjecture, we believe the following holds.

\noindent
{\bf Conjecture}  Let $k \geq 5$.  Then
$$
RR(x+y+kz = 3w) = 
\left(\left\lfloor
\frac{k+4}{3}\right\rfloor\right)^2 
+\left\{
\begin{array}{rl}
-\frac{k}{9}&\mbox{if } k \equiv 0 \,(\mbox{mod }9)\\
0&\mbox{if } k \equiv 1,6 \,(\mbox{mod }9)\\
-\frac{k+7}{9}-1&\mbox{if } k \equiv 2 \,(\mbox{mod }9)\\
-1&\mbox{if } k \equiv 3, 8 \,(\mbox{mod }9)\\
1&\mbox{if } k \equiv 4\,(\mbox{mod }9)\\
-\frac{k+4}{9}&\mbox{if } k \equiv 5 \,(\mbox{mod }9)\\
\frac{k+2}{9}&\mbox{if } k \equiv 7 \,(\mbox{mod }9)
\end{array}.
\right.
$$

\section*{\normalsize Acknowledgment}  We thank Dan Saracino for
a very careful reading that caught some errors in
a previous draft.

\vskip 30pt

We end with a table of calculated values of $RR(\mathcal{E}(k,j))$ for
small values of $k$ and $j$.  These were calculated by a standard
backtrack algorithm.  We thank Joey Parrish for helping with
implementation efficiency of the algorithm.  The program can be
downloaded as {\tt RADONUMBERS} at the second author's
homepage ({\tt
http://math.colgate.edu/$\sim$aaron/programs.html}).

\tiny
$$
\begin{array}{r|rrrrrrrrrrrrrrrrrrrrrrrr}
&\mathbf{\ell=1}&\mathbf{2}&\mathbf{3}&\mathbf{4}&\mathbf{5}&
\mathbf{6}&\mathbf{7}&\mathbf{8}&\mathbf{9}&\mathbf{10}&
\mathbf{11}&\mathbf{12}&\mathbf{13}&\mathbf{14}&\mathbf{15}&
\mathbf{16}&\mathbf{17}&\mathbf{18}&\mathbf{19}&\mathbf{20}
&\mathbf{21}
&\mathbf{22}
&\mathbf{23}\\ \hline

k\mathbf{=1}&
11&4&1&4&9&
4&10&12&14&16&
18&20&22&24&26&
38&40&43&48&50&53&59&62&
\\

\mathbf{2}&
19&5&4&1&4&
3&4&5&7&8&
9&9&17&18&20&
21&23&24&26&27&32&33&35&
\\

\mathbf{3}&
29&8&5&4&1&
4&3&4&9&6&
7&10&9&10&9&
13&14&15&16&16&26&27&28&
\\

\mathbf{4}&
41&9&4&4&5&
1&4&3&6&5&
6&4&10&8&11&
10&11&9&13&14&14&15&16&
\\

\mathbf{5}&
55&15&8&6&8&
5&1&4&6&8&
5&6&4&10&11&
11&11&13&12&14&10&14&15&
\\

\mathbf{6}&
71&17&9&5&6&
6&5&1&9&5&
8&5&11&4&9&
10&11&15&13&12&16&10&12&
\\

\mathbf{7}&
89&23&10&7&6&
7&11&5&1&4&
5&8&6&11&7&
10&12&12&13&13&14&15&10&
\\

\mathbf{8}&
109&25&15&8&6&
9&8&8&5&1&
6&5&8&5&10&
7&10&9&11&12&16&17&15&
\\

\mathbf{9}&
131&34&15&9&9&
9&9&7&14&5&
1&9&6&8&10&
10&12&12&9&11&15&12&14&
\\

\mathbf{10}&
155&37&16&10&10&
9&8&7&8&10&
5&1&9&5&8&
10&10&9&8&9&12&11&12&
\\

\mathbf{11}&
181&46&22&15&10&
11&8&11&9&10&
17&5&1&6&5&
8&6&11&9&8&11&16&11&
\\

\mathbf{12}&
209&49&24&16&11&
9&10&10&12&8&
10&12&5&1&9&
6&8&9&11&9&12&10&12&
\\

\mathbf{13}&
239&61&26&17&10&
11&13&13&10&10&
13&10&20&5&1&
9&6&6&9&11&7&10&13&
\\

\mathbf{14}&
271&65&34&18&14&
11&14&11&12&13&
8&9&11&14&5&
1&9&5&6&9&14&7&9&
\\

\mathbf{15}&
305&77&36&24&15&
13&13&15&14&15&
12&14&12&12&23&
5&1&9&6&12&9&10&7&
\\

\mathbf{16}&
341&81&38&24&16&
11&14&12&14&13&
12&12&11&10&12&
16&5&1&9&10&8&9&12&
\\

\mathbf{17}&
379&96&45&24&17&
14&15&15&17&12&
19&15&14&13&17&
14&26&5&1&9&8&8&9&
\\

\mathbf{18}&
419&101&47&25&18&
15&14&15&18&15&
12&12&16&10&12&
11&15&18&5&1&9&6&8&
\\

\mathbf{19}&
461&116&49&32&24&
17&14&17&18&18&
16&20&15&14&15&
13&15&15&29&5&1&9&8&
\\

\mathbf{20}&
505&121&60&34&25&
15&16&16&18&15&
14&20&17&19&15&
14&14&12&15&20&5&1&9&
\\

\mathbf{21}&
551&139&63&36&26&
18&19&16&15&20&
20&17&19&18&18&
14&14&17&21&16&32&5&1&
\\

\mathbf{22}&
599&145&65&37&27&
18&17&15&17&19&
18&19&15&17&15&
17&16&14&14&13&18&22&5&
\\

\mathbf{23}&
649&163&78&45&28&
23&18&18&19&19&
23&22&17&22&17&
21&14&18&16&16&17&18&35&
\\

\mathbf{24}&
701&169&\geq 81&47&32&
24&19&16&17&20&
21&18&20&16&20&
20&18&16&14&18&15&14&19&
\\

\mathbf{25}&
755&190&\geq 84&49&35&
25&19&20&16&20&
23&21&25&16&28&
23&19&18&21&18&16&17&25&
\\

\mathbf{26}&
811&197&\geq 94&51&35&
26&17&19&20&20&
23&21&22&21&17&
18&20&18&19&15&17&17&17&
\\

\mathbf{27}&
869&218&\geq 97&62&35&
27&23&22&27&22&
21&20&23&26&20&
21&22&27&19&17&21&19&19&
\\

\mathbf{28}&
929&225&\geq 100&64&36&
27&25&21&21&18&
23&22&22&21&16&
18&20&21&19&21&21&16&16&
\\

\mathbf{29}&
991&249&\geq 115&66&44&
35&26&22&21&23&
22&24&27&25&26&
19&30&21&23&25&17&21&19&
\\

\mathbf{30}&
1055&257&\geq 118&68&46&
36&24&19&23&22&
22&21&25&24&30&
17&20&24&22&21&22&18&16&
\\

\mathbf{31}&
1121&281&\geq 122&75&47&
37&26&22&23&24&
22&25&25&26&22&
28&20&30&22&20&19&24&21&
\\

\mathbf{32}&
1189&289&\geq 138&77&49&
38&27&24&24&23&
25&24&26&24&27&
24&18&20&21&22&24&25&19&
\\

\mathbf{33}&
1259&316&\geq 141&\geq 79&50&
39&26&25&25&23&
27&26&25&25&24&
24&29&20&29&29&27&28&18&
\\

\mathbf{34}&
1331&325&\geq 147&\geq 81&60&
40&34&24&23&24&
26&23&27&28&29&
27&28&20&21&26&23&26&24&
\\

\mathbf{35}&
1405&352&\geq 161&\geq 93&62&
45&35&26&25&26&
26&27&27&28&28&
26&25&34&21&33&30&25&26&
\\

\mathbf{36}&
1481&361&\geq 165&\geq 96&63&
45&36&24&27&26&
27&30&27&24&30&
27&29&27&20&22&27&22&27&
\\

\mathbf{37}&
1559&391&\geq 169&\geq 98&65&
49&37&26&26&26&
30&28&25&27&28&
31&27&27&36&22&31&26&25&
\\

\mathbf{38}&
1639&401&\geq 188&\geq 100&65&
48&38&28&28&25&
26&27&30&30&29&
30&31&32&38&22&23&28&24&
\\

\mathbf{39}&
1721&431&\geq 193&\geq 113&77&
48&39&35&29&28&
28&31&30&29&32&
32&31&29&28&32&23&31&31&
\\

\mathbf{40}&
1805	&441&\geq 196&\geq 116&79&
49&39&35&29&25&
28&26&31&27&31&
32&31&32&32&30&22&24&29&

\end{array}
$$
\vskip -10pt
$$
\begin{array}{r|rrrrrrrrrrrrrrrrrrrrrr}
&
\mathbf{24}&\mathbf{25}&
\mathbf{26}&\mathbf{27}&\mathbf{28}&\mathbf{29}&\mathbf{30}&
\mathbf{31}&\mathbf{32}&\mathbf{33}&\mathbf{34}&\mathbf{35}&
\mathbf{36}&\mathbf{37}&\mathbf{38}&\mathbf{39}&\mathbf{40}\\ \hline

k\mathbf{=1}&
64&84&87&\geq 91&\geq 98&\geq 102&\geq 105&\geq 114&\geq 118&
\geq 121&\geq 148&\geq 152&\geq 157&\geq 167&\geq 171&\geq 176&\geq 187
\\

\mathbf{2}&
36&44&
46&48&49&66&68&
\geq 70&\geq 73&\geq 77&\geq 80&\geq 82&\geq 84&\geq 93&\geq 95&\geq 98&
\geq 100
\\

\mathbf{3}&
30&31&
33&34&35&37&38&
44&45&47&48&49&
58&60&61&63&64
\\

\mathbf{4}&
17&21&
22&23&24&25&25&
37&38&39&40&42&
43&45&46&47&48
\\

\mathbf{5}&
14&17&
15&16&16&21&22&
23&23&24&25&25&
31&32&33&34&35
\\

\mathbf{6}&
9&16&
13&18&14&15&15&
16&17&21&22&22&
23&24&24&25&25
\\

\mathbf{7}&
18&9&
16&11&18&19&14&
21&22&16&17&18&
19&21&22&22&23
\\

\mathbf{8}&
16&18&
9&19&12&17&18&
19&20&21&15&23&
24&18&19&24&16
\\

\mathbf{9}&
16&16&
18&27&19&16&21&
18&19&20&21&22&
26&24&16&26&27
\\

\mathbf{10}&
13&18&
17&16&17&19&20&
21&17&23&19&20&
21&22&23&24&25
\\

\mathbf{11}&
12&15&
12&15&16&17&18&
20&20&21&21&19&
25&21&22&23&24
\\

\mathbf{12}&
11&14&
13&12&12&20&15&
16&16&20&21&20&
21&22&20&26&22
\\

\mathbf{13}&
12&10&
20&13&12&12&13&
15&16&16&17&18&
20&20&21&22&22
\\

\mathbf{14}&
10&10&
10&12&13&15&12&
13&17&21&16&20&
18&19&19&20&21
\\

\mathbf{15}&
10&15&
13&10&12&14&23&
15&14&15&15&28&
20&17&19&21&22
\\

\mathbf{16}&
7&11&
9&12&11&14&12&
16&15&18&14&16&
14&23&18&20&20
\\

\mathbf{17}&
12&8&
13&9&11&12&14&
12&18&14&26&17&
16&16&17&23&20
\\

\mathbf{18}&
9&10&
7&27&9&16&11&
18&13&15&14&16&
17&29&12&18&17
\\

\mathbf{19}&
10&10&
10&9&13&13&16&
12&17&13&15&14&
16&16&29&19&18
\\

\mathbf{20}&
10&15&
9&10&8&13&14&
16&12&13&15&20&
16&16&16&20&19
\\

\mathbf{21}&
9&10&
9&9&16&12&15&
11&17&13&17&18&
18&15&16&16&20
\\

\mathbf{22}&
1&9&
8&10&7&10&12&
13&11&20&13&12&
13&19&15&25&16
\\

\mathbf{23}&
5&1&
9&9&9&7&10&
12&13&11&16&17&
13&13&17&15&17
\\

\mathbf{24}&
24&5&
1&9&9&10&7&
15&14&13&11&14&
18&14&15&20&17
\\

\mathbf{25}&
19&38&
5&1&9&10&15&
10&12&14&13&15&
14&18&20&13&20
\\

\mathbf{26}&
15&20&
26&5&1&9&10&
10&10&12&11&14&
11&15&17&26&13
\\

\mathbf{27}&
18&21&
22&41&5&1&9&
10&10&15&11&13&
27&14&16&17&18
\\

\mathbf{28}&
16&19&
16&22&28&5&1&
9&10&10&10&21&
12&18&13&16&17
\\

\mathbf{29}&
18&19&
20&29&22&44&5&
1&9&10&10&12&
12&11&16&13&17
\\

\mathbf{30}&
21&20&
20&20&17&23&30&
5&1&9&10&11&
15&12&11&18&14
\\

\mathbf{31}&
25&22&
19&21&21&23&24&
47&5&1&9&10&
11&10&11&11&16
\\

\mathbf{32}&
20&26&
18&17&20&20&18&
24&32&5&1&9&
10&11&13&12&16
\\

\mathbf{33}&
20&23&
21&20&25&22&21&
33&26&50&5&1&
9&10&10&14&8
\\

\mathbf{34}&
19&23&
18&22&18&21&18&
21&19&27&34&5&
1&9&10&11&11
\\

\mathbf{35}&
25&22&
22&27&25&23&23&
23&25&25&27&53&
5&1&9&10&12
\\

\mathbf{36}&
24&23&
21&27&24&23&21&
21&22&23&20&27&
36&5&1&9&10
\\

\mathbf{37}&
28&23&
22&24&29&25&23&
23&22&24&25&37&
28&56&5&1&9
\\

\mathbf{38}&
27&24&
26&22&26&24&23&
24&21&23&22&25&
21&30&38&5&1
\\

\mathbf{39}&
28&29&
33&25&23&22&28&
27&25&26&25&25&
27&29&30&59&5
\\

\mathbf{40}&
25&28&
24&22&23&28&25&
26&24&27&22&28&
22&27&22&31&40

\end{array}
$$
\centerline{\bf Table 1:  Some Values of $RR(x+y+kz = \ell w)$}

\normalsize

\section*{\normalsize References}
\footnotesize
\parindent=0pt

[BB] A. Beutelspacher and W. Brestovansky,
Generalized Schur numbers, Combinatorial theory (Schloss
Rauischholzhausen, 1982), pp. 30--38,  Lecture Notes in Math., 969,
Springer, Berlin-New York, 1982.

[BL] S. Burr and S. Loo, On Rado Numbers II, unpublished.

[GS] S. Guo and Z-W. Sun,
Determination of the two-color Rado number for $a_1x_1+\cdots+a_mx_m=x_0$,
to appear in {\it JCTA}, preprint available at
{\tt arXiv://math.CO/0601409}.

[HS] B. Hopkins and D. Schaal,
On Rado Numbers for $\sum_{i=1}^{m-1}a_ix_i=x_m$,
{\it Adv. Applied Math.} {\bf 35} (2005), 433-441.

[MR] K. Myers and A. Robertson, Two Color Off-diagonal Rado-type
Numbers, preprint available at \vskip -8pt
{\tt arXiv://math.CO/0606451}.

[R] R. Rado, Studien zur Kombinatorik, {\it Mathematische Zeitschrift}
{\bf 36} (1933), 424-480.

\end{document}